\documentstyle[12pt, amssymb]{article}
\begin{document}
\newcommand{\bref}[1]{$\mbox{(\ref{#1})}$}
\newcommand{\g}{\gamma}
\newcommand{\s}{\sum\limits_{k =1}^{d}}
\newcommand{\r}{{\Bbb {R}} ^d}
\newcommand{\ec}{\left(\frac{|E|}{2C}\right)}

\begin{center}
\bf Some results related to the Logvinenko-Sereda theorem

~

\rm Oleg Kovrijkine\\
Department of Mathematics  253-37\\
California Institute of Technology\\
Pasadena, CA 9125, USA\\
{\sl E-mail address}: olegk@its.caltech.edu\\
\end{center}

\begin{abstract}We prove several results related to the  theorem of Logvinenko
and Sereda on  determining sets for functions with Fourier transforms supported
 in an interval. We obtain a polynomial instead of exponential bound in this
theorem, and we extend it to the case of functions with Fourier transforms
supported in the union of a bounded number of intervals.
The same results hold in all dimensions.
\end{abstract}

	The purpose of this work is to study the behavior of 
functions whose Fourier transforms are supported in an interval
 or in a union of finitely many intervals on \lq\lq thick" subsets 
of the real line. A main result of this type was proven by Logvinenko and Sereda.\\

By a \lq\lq thick" subset of the real line we mean a measurable
 set $E$ for
 which there exist $a > 0$ and $\gamma > 0$ such
 that
\begin{eqnarray}|E \cap I| \ge \gamma \cdot a \label{e}
\end{eqnarray}
for every interval $I$ of length $a$.

~

{\bf The Logvinenko-Sereda Theorem}:
{\it let $J$ be an interval with $|J| = b$. If $f \in L^p$, $p \in [1, + \infty]$,
 and
 supp $\hat f \subset J$ and if a measurable set $E$ satisfies \bref {e}
 then 
\begin{eqnarray}\|f\|_{L^p (E)} \ge \exp(-C \cdot \frac
 {(ab + 1)} {\g})\cdot\|f\|_p.\label{ls}\end{eqnarray} }

It is a well-known fact that the condition \bref{e} is also necessary for 
an inequality of the form
$$\|f\|_{L^p (E)} \ge C\cdot\|f\|_p$$
to hold. See for example (\cite{HJ}, p.113).\\

We will improve the estimate \bref{ls} by getting a polynomial
 dependence on $\gamma$ and show that our estimate is optimal except for the
 constant $C$:\\
{\bf Theorem 1}:
$$\|f\|_{L^p (E)} \ge {\left
 ( \frac {\gamma} {C} \right )}^{C \cdot (ab + 1)}\cdot\|f\|_p. $$

	We will also generalize the Logvinenko-Sereda theorem to functions
 whose Fourier transforms are supported on a union of finitely many
 intervals:\\

{\bf Theorem 2}:
{\it let $J_k$ be intervals with $|J_k| = b$. If $f \in L^p$, $p \in
 [1, + \infty]$,
 and supp $\hat f \subset \bigcup \limits_1^n J_k$ and if a measurable
 set $E$
 satisfies \bref {e} then 
$$\|f\|_{L^p (E)} \ge c(\gamma, n, ab, p)\cdot\|f\|_p$$
where $c(\gamma, n, ab, p) = {\left ( \frac C {\gamma} \right )}^{-ab
\left ( \frac C {\gamma} \right )^n -n + \frac {p-1} p}$ depends only on the number
 of intervals but not how they are placed.}

~

	Note that the constant $C$ below is not fixed and varies appropriately
 from one equality or inequality to another even without mentioning it.

~

Proof of {\bf Theorem 1}:\\
	First we treat the case when $p \in [1, + \infty)$.
	Without loss of generality we can always assume that
 $J = [-\frac b 2, \frac b 2]$.
By considering $f(\frac x a)$ instead of $f$ we can also assume
 that $|E \cap [x,x+1]| \ge \gamma$ $\forall x$ and
 supp $\hat f \subset [-\frac {ab} 2,\frac {ab} 2]$, just say
  supp $\hat f \subset [-\frac b 2, \frac b 2]$. Bernstein's inequality
 (\cite{B}, Theorem 11.3.3) gives that 
 \[\int |f^{(\alpha)}|^p \le (C\cdot b)^{\alpha p} \cdot \int |f|^p\]
with $C = \frac 12$.\\
Divide the whole $\Bbb R$ into intervals of length 1. Choose $A > 1$.
 Call an
 interval I bad if $\exists \alpha \ge 1$ such that 
\[\int \limits_I |f^{(\alpha)}|^p \ge A^{\alpha p}(C\cdot b)^{\alpha p}
 \cdot 
\int
 \limits_I |f|^p.\]
Then 
\begin{eqnarray} \int \limits_{\bigcup \limits_{I\; is \; bad}I}|f|^p
 &\le& 
\int 
\limits_{\bigcup \limits_{I\; is \; bad}I} \sum 
\limits_{\alpha =1}^{\infty} 
\frac 1{A^{\alpha p}(C\cdot b)^{\alpha p}}|f^{(\alpha)}|^p\nonumber
\\
 &=& \sum \limits_{\alpha =1}^{\infty} \frac 1{A^{\alpha p}
(C\cdot b)^{\alpha p}}
 \int
 \limits_{\bigcup \limits_{I\; is \; bad}I} |f^{(\alpha)}|^p\nonumber \\
&\le& \sum \limits_{\alpha =1}^{\infty} \frac 1{A^{\alpha p}
(C\cdot b)^{\alpha p}} 
\int
 |f^{(\alpha)}|^p\nonumber\\
&\le& \sum \limits_{\alpha =1}^{\infty} \frac 1{A^{\alpha p}} 
\int |f|^p\nonumber
\\
 &=& \frac 1{A^p - 1}\int |f|^p.\label{xyz2}\end{eqnarray}
Choose $A = 3$ and apply \bref{xyz2}. So 
$$\int \limits_{\bigcup \limits_{I\; is \; bad}I}|f|^p \le \frac 1 {2}
\int |f|^p.$$
Therefore 
\begin{eqnarray}\int \limits_{\bigcup \limits_{I \; is \; good}I}|f|^p \ge 
\frac 1 {2}\int |f|^p.\label{good}\end{eqnarray}
 We claim that $\exists B > 1$ such that if $I$ is a good interval then
 $\exists x \in I$ with the property that  
$$|f^{(\alpha)}(x)|^p \le 2\cdot B^{\alpha p} (C\cdot b)^{\alpha p}
\cdot \int \limits_I |f|^p \; \; \; \forall \alpha \ge 0.$$
 Suppose towards a contradiction that this is not true. Then  

\begin{equation}2\cdot \int \limits_I |f|^p \le \sum \limits_{\alpha =0}^{\infty}
\frac 1 {B^{\alpha p}(C\cdot b)^{\alpha p}} |f^{(\alpha)}(x)|^p \; \; \; \forall x
 \in I. \label{xyz}\end{equation}
Integrate both sides of \bref{xyz} over $I$:
\begin{eqnarray}2 \cdot \int \limits_I |f|^p &\le& \sum \limits_{\alpha =0}^{\infty}
 \frac 1 {B^{\alpha p}(C\cdot b)^{\alpha p}} \int \limits_I |f^{(\alpha)}(x)|^p\nonumber
 \\
 &\le& \sum \limits_{\alpha =0}^{\infty} \frac 1 {B^{\alpha p}} \int \limits_I |f|^p\nonumber
 \\
 &=& \frac 1 {1 - {\left ( \frac 1 {B} \right )}^p}\int \limits_I |f|^p.
\label{xyz3}\end{eqnarray}
Choose $B = 3$ and apply \bref {xyz3}. So 
$$2 \cdot \int \limits_I |f|^p \le \frac 3 2 \int \limits_I |f|^p.$$
This contradiction proves our claim.\\
We will need to prove the following local estimate:
$$\int\limits_{E \cap I} |f|^p \ge \left ( \frac {\gamma} C \right )^ {Cbp + 2}
 \int\limits_{ I} |f|^p$$
for every good interval $I$.
	Without loss of generality we can assume that $I = [-\frac 1 2,\frac 1 2]$
 by considering a shift $f(x-n)$ which has supp$\widehat {f(x-n)} \subset
 [-\frac b 2,\frac b 2]$. Therefore if $y \in D(0,R) \subset D(x,R+\frac 1 2)$ then 
 \begin{eqnarray}|f(y)| &\le& \sum \limits_{\alpha =0}^{\infty}
 \frac {|f^{(\alpha)}(x)|} {\alpha !} \cdot |y - x| ^ {\alpha}\nonumber
 \\
&\le& \sum \limits_{\alpha =0}^{\infty}2^{\frac 1 p} \frac {(R+\frac 1 2)^{\alpha}
\cdot (Cb)^{\alpha}} {\alpha !}\|f\|_{L^p(I)}\nonumber
\\
&=& 2^{\frac 1 p}\exp(Cb(R+\frac 1 2)) \cdot \|f\|_{L^p(I)}. \label{exp}
\end{eqnarray}

Now we will give a local estimate for analytic functions.\\

{\bf Lemma 1}:
{\it Let $\phi(z)$ be analytic in $D(0,5)$ and let $I$ be an interval of 
length $1$ such that $0 \in I$ and let $E \subset I$ be a measurable set
 of positive measure. If $|\phi(0)| \ge 1 $ and $M = \max \limits _{|z| \le 4}
 |\phi(z)|$ then

\begin{eqnarray}\sup \limits _{x \in I}|\phi(x)| \le \left (\frac C {|E|}
 \right)^{\frac {\ln M} {\ln 2}}\sup\limits _{x \in E}|\phi(x)|.
 \label {phi}\end{eqnarray}}

~

Proof of {\bf Lemma 1}:\\
Let $a_1,$ $a_2,$...$a_n$ be the zeros of $\phi$ in $D(0,2)$. If 
$$g(z) = \phi(z) \cdot \prod \limits _{k=1}^{n} \frac {4 - \bar a_k z} 
{2(a_k - z)}
 = \phi(z) \cdot \frac {Q(z)} {P(z)}$$
then $|g(0)| \ge 1$ and $\max \limits _{|z| \le 2} |g(z)| \le M$ by the
 property of Blaschke products. Applying Harnack's inequality 
to the
 positive harmonic function $\ln M - \ln |g(z)|$ in $D(0,2)$ we
 have:
$$\max \limits _{|z| \le 1} ( \ln M - \ln |g(z)|) \le 3
 \ln M.$$
Therefore
$$\min \limits _{|z| \le 1} |g(z)| \ge M^{-2}.$$
This gives us
$$\frac {\max \limits _{x \in I}|g(x)|} {\min \limits _{x \in I}|g(x)|} \le M^3.$$  

We can give a similar estimate for $Q$:
\begin{eqnarray*} \frac {\max \limits _{x \in I}|Q(x)|} 
{\min \limits _{x \in I}|Q(x)|}  
&\le& \frac {\max \limits _{|z| \le 1}\prod \limits _{k=1} ^{n} |4 - \bar a_k z|} 
{\min \limits _{|z| \le 1}\prod \limits_{k=1} ^{n} |4 - \bar a_k z|} \\
&\le& 3^n.\end{eqnarray*}

From the Remez inequality for polynomials (\cite{BE}, Theorem 5.1.1) it follows 
that:
$$\sup \limits _{x \in I} |P(x)| \le  \left ( \frac 4 {|E|} \right )^n
 \cdot \sup \limits _{x \in E} |P(x)|.$$
Therefore
\begin{eqnarray*} \sup \limits _{x \in I}|\phi(x)| &\le& \max \limits _{x \in I}
|g(x)| \cdot \frac {\max \limits _{x \in I} |P(x)|} {\min \limits _{x \in I}|Q(x)|} 
\\
&\le& M^3 \cdot 3^n \cdot \left ( \frac C {|E|} \right ) ^n \cdot \min \limits _{x \in I}
|g(x)| \cdot \frac {\sup \limits _{x \in E} |P(x)|} {\max \limits _{x \in I}|Q(x)|} 
\\
&\le& M^3 \cdot \left ( \frac C {|E|} \right ) ^n \cdot\sup \limits _{x \in E}|\phi(x)|.\end{eqnarray*}
From Jensen's formula it follows that $n \le \frac {\ln M} {\ln 2}.$
Therefore
$$\sup \limits _{x \in I}|\phi(x)| \le \left (\frac C {|E|} \right)^{\frac {\ln M} 
{\ln 2}}\sup\limits _{x \in E}|\phi(x)|.$$
\hfill$\square$

~

{\bf Corollary}: if $p \in [1, \infty)$ then
\begin{eqnarray}\| \phi \|_{L^{p}(I)} \le \left (\frac C {|E|} \right)^{\frac {\ln M} 
{\ln 2} + \frac 1 p}\| \phi\|_{L^{p}(E)}. \label {phi1}\end{eqnarray}

It follows from \bref {phi} that:
$$| \{x \in I : |\phi(x)| < \left (\frac  {\epsilon} C \right)^{\frac {\ln M} {\ln 2}}\|\phi\|_{L^{\infty}(I)}\}| \le \epsilon \; \; \; \epsilon > 0.$$
If we put $\epsilon = \frac {|E|} 2$ then 
$$| \{x \in I : |\phi(x)| < \left (\frac  {|E|} {2C} \right)^{\frac {\ln M} {\ln 2}}\|\phi\|_{L^{\infty}(I)}\}| \le \frac {|E|} 2.$$
Therefore
\begin{eqnarray*} \int \limits _E |\phi|^p &\ge& \int \limits _E \chi_{|\phi|
 \ge \ec^{\frac {\ln M} {\ln 2}}\cdot \|\phi\|_{L^{\infty}(I)}} \cdot |\phi|^p \\
&\ge& \frac {|E|} 2 \cdot \left ( \frac {|E|} {2C}  \right )^ {p\frac {\ln M} {\ln 2}}
\cdot \|\phi\|_{L^{\infty}(I)}^p \\
&\ge& \left ( \frac {|E|} {2C}  \right )^ {p\frac {\ln M} {\ln 2} + 1} \cdot \int 
\limits _I |\phi|^p.
\end{eqnarray*}\hfill$\square$

~

Now we are in a position to proceed with the proof of our theorem. We can assume that 
$\int _I |f|^p = 1$. Then $\exists x_0 \in I$ such that
$|f(x_0)| \ge 1$.
 Applying \bref {phi1} to $\phi (z) = f (z + x_0)$, $I - x_0$ and $(E \cap I) -x_0$ 
we have:
$$\int _{E \cap I} |f|^p \ge \left ( \frac {|E \cap I|} C \right )^ {p\frac {\ln M} 
{\ln 2} + 1}\int _{I} |f|^p.$$
Apply \bref {exp} to get
\begin{eqnarray*} M &\le&  \max \limits _{|z| \le 4 + \frac 12} |f(z)| \\
&\le& 2^{\frac 1 p}\exp(5Cb).\end{eqnarray*}
Therefore  we have :
\begin{eqnarray}
\int _{E \cap I} |f|^p \ge \left ( \frac {\gamma} C \right )^ {Cbp + 2} \int _{ I} 
|f|^p.
 \label {loc}\end{eqnarray}
 Summing \bref {loc} over all good intervals and applying \bref{good} we have
\begin{eqnarray*}\int \limits _E |f|^p &\ge& \int \limits _{E \cap \bigcup 
\limits_{I \; is \; good} I} |f|^p \\
&\ge& \left ( \frac {\gamma}{C}\right ) ^{Cbp + 2}\cdot \int \limits _{\bigcup 
\limits_{I \; is \; good} I} |f|^p \\
&\ge& \frac12 \left ( \frac {\gamma}{C}\right ) ^{Cbp + 2}\cdot \int |f|^p.
\end{eqnarray*}
Replacing $b$ with $ab$ and choosing a new $C$ we have:
$$\int\limits_{E} |f|^p \ge {\left ( \frac {\gamma} {C} \right )}^{Cabp  +2}\cdot
\int |f|^p.$$

If $p = \infty$ then the proof is almost the same: $\|f\|_{L^{\infty} (\bigcup 
\limits_{I \; is \; good} I)} = \|f\|_{\infty}$. If $I$ is good then $\|f\|_{L^{\infty} 
(E \cap I)} 
\ge \left ( \frac {\gamma}{C}\right ) ^{Cb + 1}\cdot \|f\|_{L^{\infty} (I)}$. 
Hence
 $$\|f\|_{L^{\infty} (E)} \ge \left ( \frac {\gamma}{C}\right ) ^{Cb + 1}\cdot 
\|f\|_{\infty}.$$
End of proof of {\bf Theorem 1}.

~

If we keep track of all the constants and do the calculations more accurately then 
we can get
 that if $p \in [1, \infty)$:
$$\|f\|_{L^p (E)} \ge {\left ( \frac {\gamma} {300} \right )}^{33ab + \frac 2 p}
\cdot\|f\|_p,$$
if $p = \infty$:
$$\|f\|_{L^{\infty} (E)} \ge {\left ( \frac {\gamma} {100} \right )}^{33ab + 1 }\cdot\|f\|_{\infty}.$$
However, if we try to minimize the factor in front of $ab$ then we can get the 
following 
estimate:
$$\|f\|_{L^p (E)} \ge {\left ( \frac {\gamma} {C} \right )}^{(\frac {(1 + e)} 2 + 
\epsilon)
\cdot ab + A(\epsilon)}\cdot\|f\|_p \;\;\; \forall\epsilon > 0.$$

The following example suggests that the right behavior of the estimate in the 
Logvinenko-Sereda Theorem is $\gamma$ to the power of a linear function of $ab$ 
and how far we are from the exact factor in front of $ab$:\\
 Let $E$ be a 1-periodic set such that 
$$E \cap [-\frac12, \frac12] = [-\frac12, -\frac12 + \frac {\gamma}2] \cap  
[\frac12 - \frac {\gamma}2, \frac12]$$
and let 
$$f(x) = \left ( \frac {\sin (2\pi x)} x \right ) ^ {[\frac b {4\pi}]}$$
If $b$ is large enough we have:
$$\|f\|_{L^p(E)} \le {\left ( \frac {\gamma} {C} \right )}^{\frac b {4\pi} -1}\|f\|_p$$
and supp$\hat f \subset [-\frac b 2,\frac b 2]$.

~

{\bf Remark 1}: when $ab$ is sufficiently small the proof of the theorem is much simpler:
 if $ab \le 1$ then $\|f\|_{L^p (E)} \ge \frac {\gamma^{\frac1p}} 2\|f\|_p$. This can
 be proven very easily. If $p \in [1, + \infty)$ we have \begin{eqnarray*}|f(x)|^p = 
|f(y) - \int_x^y f'(t)dt|^p  &\ge& \frac {|f(y)|^p}{2^{p-1}} - |\int_x^y f'(t)dt|^p 
\\
 &\ge& \frac {|f(y)|^p}{2^{p-1}} - \int \limits_I |f'|^p \cdot a^{p-1}
\end{eqnarray*}
 where $x,y \in I$, $|I| = a$. Hence 
\begin{eqnarray*}a\cdot\int \limits _{E\cap I}|f(x)|^p dx &=& \int\limits_I
 (\int 
\limits _{E\cap I}|f(x)|^p dx)dy \\
&\ge& |{E\cap I}| \cdot \left ( \frac 1 {2^{p-1}} \int \limits _ I |f|^p - a^p
\int 
\limits_I |f'|^p \right ).\end{eqnarray*}
Therefore $\frac 1 {\gamma} \cdot\int \limits _{E\cap I}|f|^p \ge \frac 1 {2^{p-1}} 
\int 
\limits _ I |f|^p - a^p\int \limits_I |f'|^p $.
Summing over all intervals $I$ we have 
\begin{eqnarray*}\frac 1 {\gamma}\cdot\int \limits _E |f|^p &\ge& \frac 1 {2^{p-1}} 
\int |f|^p - a^p\int |f'|^p \\
&\ge& \frac 1 {2^{p-1}} \int |f|^p - (\frac b 2)^p a^p \int |f|^p \\
&\ge& \frac 1 {2^p}\cdot \int |f|^p.\end{eqnarray*}
Using $\|(f^p)'\|_1 \le \frac {pb} 2 \|f^p\|_1$ we can get a similar result. The proof 
for $p = \infty$ is even easier. \\
In a similar way  we can treat the case when $1 - \gamma$ is sufficiently small 
depending on
 $ab$: if $p \in [1, \infty)$ and $1 - \gamma \le \frac 1 {2 + pab}$ then 
$\|f\|_{L^p (E)}^p \ge \frac 12 \|f\|_p^p$.

~

Proof of {\bf Theorem 2}:\\
      Let $J_k = [\lambda_k - \frac b 2, \lambda_k + \frac b 2]$,  $k = 1, 2, ..., n$.
	First we will prove a special case of {\bf Theorem 2}:

~

 {\bf Theorem 2$'$}:
{\it  
if $\lambda_{k+1} - \lambda_{k} \ge 2b > 0$   $(k = 1,2, ...,n - 1)$\\
then 
$$\|f\|_{L^p(E)} \ge c(\gamma, n, ab, p) \cdot \|f\|_{L^p}$$
where $c(\gamma, n, ab, p) = \left (\frac {\gamma} C \right)^ {ab 
\left ( \frac C {\gamma}
 \right )^n + n - \frac {p-1}p}$.}

~

Proof of {\bf Theorem 2$'$}:\\
	Let $\hat f (x) = \sum \limits _{k = 1} ^n \hat f_k (x - \lambda_k)$ where 
supp$\hat f_k \subset [- \frac b 2,\frac b 2]$ and $f (x) = \sum \limits _{k = 1} ^n 
f_k (x) e^{i\lambda_k x}$. 
	The following lemma gives an estimate of $\|f_k\|_p$ from above.

~

{\bf Lemma 2}:
{\it \begin{eqnarray}\|f_k\|_p \le C\|f\|_p   \;\;\;(k = 1,2,...,n)  .
\label{fk}\end{eqnarray}}
 
Proof of {\bf Lemma 2}:\\
Let $\phi$ be a Schwartz function such that supp $\hat \phi \subset [-1,1]$ and
 $\hat \phi (x) = 1$ for $x \in [-\frac 1 2,\frac 1 2]$. Then $\hat f_k (x) =
 \hat f \cdot \hat \phi (\frac {x - \lambda_k} {b})$. Therefore $f_k =
 f*(b\phi(bx)e^{i\lambda_k x})$. Applying Young's inequality we have
 $\|f_k\|_p \le \|f\|_p \cdot \|\phi\|_1$.
 \hfill$\square$

~

	We will also need the following auxilary lemma:

~

{\bf Lemma 3}:
{\it if $r(x) = \sum \limits _{k = 1} ^{n} p_k(x)e^{i\lambda_k x}$ \\
 where $p_k(x)$ is a polynomial of degree $\le m-1$ and $E \subset I$ is 
measurable with $|E| > 0$ then 
\begin{eqnarray}\|r\|_{L^p(I)} \le \left ( \frac {C|I|} {|E|} \right )^ {nm - 
\frac {(p -1)} p}\cdot \|r\|_{L^p(E)} .\label {r}\end{eqnarray}}

Proof of {\bf Lemma 3}:\\
	First we prove the statement for pure trigonometric polynomials, i.e. \\
if
$g(x) = \sum \limits _{k = 1} ^{n} c_k e^{i\lambda_k x}$ then
\begin{eqnarray}\|g\|_{L^p(I)} \le \left ( \frac {C|I|} {|E|} \right )^{n - 
\frac {(p -1)} p}
\cdot \|g\|_{L^p(E)}.\label {g}\end{eqnarray}

This follows from a theorem on trigonometric polynomials by F. Nazarov 
(\cite{N}, Theorem 1.5)
 saying that:
\begin{eqnarray}\|g\|_{L^{\infty}(I)} \le \left ( \frac {C|I|} {|E|} \right )^{n - 1}
\cdot \|g\|_{L^{\infty}(E)}.\label{naz}\end{eqnarray}
An argument similar to the proof of the {\bf Corollary} to {\bf Lemma 1} shows that 
\bref{g} follows from \bref{naz}. \\
 
If $p(x) = \sum \limits _{l = 0} ^{m-1} a_l x^l$ is a polynomial of degree 
$m-1$
 then it can be approximated uniformly on an interval with a trigonometric 
polynomial
 of order $\le m$  
$$\tilde {p}(x) = \sum \limits _{l = 0} ^{m-1} a_l \left ( \frac {e^{i\lambda x} - 1}
 {i\lambda} \right )^l = \sum \limits _{l = 0} ^{m-1} \tilde {a}_l e^{il\lambda x}$$
because $x = \lim \limits _{\lambda \rightarrow 0} 
\frac {e^{i\lambda x} - 1} {i\lambda}$
 uniformly on an interval. Applying \bref {g} to the trigonometric polynomial
 of order $mn$
$$\tilde r(x) = \sum \limits _{k = 1} ^{n} \tilde p_k(x)e^{i\lambda_k x}$$
 and taking the limit we have the desired result: 
$$\|r\|_{L^p(I)} \le \left ( \frac {C|I|} {|E|} \right )^ {nm - \frac {(p -1)} p}
\cdot \|r\|_{L^p(E)} .$$ \hfill$\square$

~

	We will need the Taylor formula:
\begin{eqnarray*} g(x) &=& \sum \limits _{l = 0} ^{m-1} \frac {g^{l} (s)} {l!}
 (x - s)^l + \frac 1 {(m-1)!}\int \limits_{s}^{x} g^{(m)} (t) (x - t)^{m-1} dt \\ 
&=& p(x) + \frac 1 {(m-1)!}\int \limits_{s}^{x} g^{(m)} (t) (x - t)^{m-1}dt 
\end{eqnarray*}
where $p(x)$ is a polynomial of degree $m-1$.

	Now we are in a position to proceed with the proof of {\bf Theorem 2$'$}.

First we assume that $p \in [1, \infty)$. Divide the whole $\Bbb R$ into intervals
 of length $a$ each. Consider one of them: $I = [s, s+a]$. Then 
\begin{eqnarray*}
f(x) &=& \sum \limits _{k = 1} ^n f_k (x) e^{i\lambda_k x}\\
&=& \sum \limits _{k = 1} ^n p_k (x) e^{i\lambda_k x} + \frac 1 {(m-1)!}
\sum \limits _{k = 1}^n e^{i\lambda_k x} \int \limits_{s}^{x} f_k^{(m)} (t) (x - t)^{m-1} dt
 \\
&=& r(x) + T(x) \end{eqnarray*}

Applying Holder's inequality we have
 \begin{eqnarray}\int \limits_I |T(x)|^p dx &\le& \frac {n^{p-1}} {[(m-1)!]^p}
 \sum \limits _{k = 1} ^n \int \limits_I |\int \limits_{s}^{x} f_k^{(m)} (t)
 (x - t)^{m-1} dt|^p dx\nonumber \\
&\le& \frac {n^{p-1}a^{pm}} {[m!]^p}\sum \limits _{k = 1} ^n \int \limits_I |f_k^{(m)}|^p.\label{t}\end{eqnarray}

\begin{eqnarray}\int \limits_I |f|^p &\le& 2^{p-1}\int \limits_I |r|^p + 2^{p-1}
\int \limits_I |T|^p\nonumber \\
 &\le& 
\left ( \frac {C|I|} {|E \cap I|} \right )^ {pnm -(p-1)}\cdot
\int \limits_{E \cap I}
 |r|^p + 2^{p-1}\int \limits_I |T|^p\nonumber \\
 &\le&
\left ( \frac {C} {\gamma} \right )^ {pnm -(p-1)}\cdot\left (2^{p-1}
\int \limits_{E \cap I} |f|^p + 2^{p-1}\int \limits_{E \cap I} |T|^p\right)
 + 2^{p-1}\int \limits_I |T|^p\nonumber \\ 
&\le& \left ( \frac {C} {\gamma} \right )^ {pnm -(p-1)}\cdot\int \limits_{E \cap I}
 |f|^p + \left ( \frac {C} {\gamma} \right )^ {pnm -(p-1)}\cdot\int \limits_{I}
 |T|^p\nonumber \\
 &\le& \left (\frac {C} {\gamma} \right )^ {pnm -(p-1)}\cdot\int \limits_{E \cap I}
 |f|^p + \left ( \frac {C} {\gamma} \right )^ {pnm -(p-1)}\cdot\frac {n^{p-1}a^{pm}}
 {[m!]^p}\sum \limits _{k = 1} ^n \int \limits_I |f_k^{(m)}|^p. \nonumber \end{eqnarray}
The second inequality is based on {\bf Lemma 3}. The last follows from \bref{t}.\\
 
Summing over all intervals $I$ we have:
\begin {eqnarray*}\int |f|^p &\le& \left (\frac {C} {\gamma} \right )^ {pnm -(p-1)}
\cdot\int \limits_{E} |f|^p + \left ( \frac {C} {\gamma} \right )^ {pnm -(p-1)}
\frac {n^{p-1}a^{pm}} {[m!]^p}\sum \limits _{k = 1} ^n \int |f_k^{(m)}|^p \\
&\le& \left (\frac {C} {\gamma} \right )^ {pnm -(p-1)}\cdot\int \limits_{E}
 |f|^p + \left ( \frac {C} {\gamma} \right )^ {pnm -(p-1)}\frac {n^{p-1}a^{pm}(Cb)^{pm}} {[m!]^p}\sum \limits _{k = 1} ^n \int |f_k|^p \\
&\le& \left (\frac {C} {\gamma} \right )^ {pnm -(p-1)}\cdot\int \limits_{E}
 |f|^p + \left ( \frac {C} {\gamma} \right )^ {pnm -(p-1)}\frac {n^p(ab)^{pm}} {[m!]^p}
 \int |f|^p\\
&\le& \left (\frac {C} {\gamma} \right )^ {pnm -(p-1)}\cdot\int \limits_{E}
 |f|^p + \left ( \frac {C} {\gamma} \right )^ {pnm}\frac {(ab)^{pm}} {m^{pm}} 
\int |f|^p.\end{eqnarray*}
The second inequality follows from Bernstein's Theorem. The third is an application
 of \bref{fk}. The last inequality is due to Stirling's formula for $m!$ and the
 fact that $n \le 2^n$.\\

Choose $m$ such that it is a positive integer and $\left (
 \frac {C} {\gamma} \right )^ {n}\frac {ab} {m} \le \frac12$, e.g. \\
 $m = 1 + [\left ( \frac {C} {\gamma} \right )^ {n}\cdot ab]$ for some $C >0$.
  Therefore 
\begin{eqnarray*}\int |f|^p &\le& \left (\frac {C} {\gamma} \right )^{pn(1 + 
\left ( \frac {C} {\gamma} \right )^ {n}\cdot ab) - (p-1)}
\cdot\int \limits_{E} |f|^p \\
&\le& \left (\frac {C} {\gamma} \right )^{ p\left ( \frac {C} {\gamma} \right )^ {n}
\cdot ab + pn - (p-1)}\cdot\int \limits_{E} |f|^p .\end{eqnarray*}

The proof for $p= \infty$ is similar and even simpler.\\
End of proof of {\bf Theorem 2$'$}.

~

	Now we can proceed with the proof of {\bf Theorem 2}. We will apply induction on
 $n$. For $n = 1$ the theorem follows from {\bf Theorem 2$'$} or the usual 
Logvinenko-Sereda Theorem. Suppose the statement is true for $n \le m $. Let $n = m + 1$.
\\ 

If $\lambda_{k+1} - \lambda_{k} \ge 2b > 0$  $(k = 1,2 ...)$ then
the result follows from {\bf Theorem 2$'$}.\\

If  $0< \lambda_{k+1} - \lambda_{k} < 2b$ for some $k$ then we can replace $b$ with
 $3b$ reducing the number of frequences $\lambda_{k}$. Therefore by induction:
\begin{eqnarray*}\|f\|_{L^p (E)} &\ge& {\left ( \frac C {\gamma} \right )}
^{-3ba{\left ( \frac C {\gamma} \right )}^m - m + \frac {p-1} p}\cdot\|f\|_p\\
 &\ge& {\left ( \frac C {\gamma} \right )}^ {-ab \left ( \frac C {\gamma} \right )
^{(m+1)} -(m+1) + \frac {p-1} p}\cdot\|f\|_p.\end{eqnarray*}
End of proof of {\bf Theorem 2}.

~

The purpose of this theorem is to prove the existence of a constant 
$c(\gamma, n, ab, p) > 0$ depending only on the number of intervals and not 
how they are placed rather than to get the best possible estimate.

~

{\bf Final remark.} By a \lq\lq thick" subset of $\r$ we mean a measurable set $E$ for
 which there exist a parallelepiped $\Pi$ with sides of length $a_1, a_2,..., a_d$ 
parallel to coordinate axes  and $\gamma > 0$ such that
\begin{eqnarray}|E \cap (\Pi + x)| \ge \gamma |\Pi| \label{eee1}\end{eqnarray}
for every $x \in \r$. {\bf Theorems 1} and {\bf 2} can be easily extended to higher
 dimensions
 with polynomial dependence on $\gamma$ for the former. The proofs are analogous to the
 previous proofs. We can assume that $\Pi$ is a unit cube. Define good cubes in a 
similar way. The main 
issue is how to obtain a local estimate for good cubes. If $|f|$ attains its maximum 
in a cube $\Pi$ at $y \in \Pi$ then following an idea of F. Nazarov we can use 
spherical coordinates 
centered at $y$ to find a segment $I$ in $\Pi$ such that $y \in I$ and 
$\frac {|E \cap I|}{|I|}
 \ge C(d) \gamma$, and reduce our problem to a 1-dimensional one. In case of 
{\bf Theorem 1} 
we can define an analytic function of one complex variable which coincides with 
$f$ on $I$. 
In case of {\bf Theorem 2} we will approximate $f$ on $I$ with a polynomial defined
 on $I$.\\

{\bf Theorem 3.} {\it let $J$ be a parallelepiped with sides of length 
$b_1, b_2,..., b_d$ 
parallel to coordinate axes. If $f \in L^p(\r)$, $p \in [1, + \infty]$, and
 supp $\hat f \subset J$ and if a measurable set $E$ satisfies \bref {eee1} 
then 
\begin{eqnarray*}\|f\|_{L^p (E)} \ge \left(\frac {\gamma}{C^d}\right)^{C(d + 
\s a_k b_k)}\|f\|_p.\end{eqnarray*} }
By an example similar to the one after {\bf Theorem 1} (with supp$\hat f$ in a 
neighborhood of a main diagonal of $J$ with the direction of ${\bf b} = (b_1,...,b_d)$
 and $E$ periodic along the same direction with period $\sim {\bf a \cdot b}/|{\bf b}|$)
 we can show that this estimate is optimal except for the constant $C$.\\

{\bf Theorem 4.} {\it let $J_l$ be parallelepipeds with sides of length 
$b_1, b_2,..., b_d$ parallel to coordinate axes. If $f \in L^p(\r)$, 
$p \in [1, + \infty]$, and
 supp $\hat f \subset \bigcup \limits_1^n J_l$ and if a measurable set 
$E$ satisfies \bref {eee1} then 
\begin{eqnarray*}\|f\|_{L^p (E)} \ge c(\gamma, n, {\bf a \cdot b}, p, d)
\|f\|_p\end{eqnarray*}
where $c(\gamma, n, {\bf a \cdot b}, p, d) = 
{\left ( \frac {C^d} {\gamma} \right )}^{-\left ( \frac {C^d} {\gamma} \right )^n 
\cdot \s a_k b_k - n + \frac {p-1} p}$ depends only on the 
number of parallelepipeds but not how they are placed. }\\

\begin{center}
{\bf Acknowledgements}
\end{center}
The author is grateful to Professor Thomas Wolff for his interest in this work.

~

\end{document}